\documentclass[12pt]{article}
\usepackage{amsmath,amssymb,amsthm,amscd,graphics,epsfig}

\newcommand{\field}[1]{\mathbb{#1}}
\newcommand{\RR}{\field{R}}

\newcommand{\al} {\alpha}

       \newcommand{\Te}{\Theta}

\newcommand{\la} {\lambda}

\newcommand{\si} {\sigma}

\newcommand{\vfi}{\varphi}
       \newcommand{\Om}{\Omega}

\newcommand{\cB}{{\mathcal B}}
\newcommand{\cF}{{\mathcal F}}
\newcommand{\cL}{{\mathcal L}}

\newtheorem{theo}{Theorem}
\newtheorem*{theoA}{Theorem A}

\newtheorem{lemm}{Lemma}[section]

\newtheorem{defi}[lemm]{Definition}
\newtheorem{prop}[lemm]{Proposition}

\newtheorem{rema}[lemm]{Remark}
\newenvironment{demo}{{\bf Proof: }}{\hfill
$\qed$\\ \medskip}

\begin{document}

\title{A Ruelle Operator for continuous time Markov
Chains }
\author{ A. Baraviera  \footnote{ Instituto de Matem\'atica, UFRGS, Porto
Alegre 91501-970, Brasil. benefici\'ario de aux\'{\i}lio financeiro
CAPES - PROCAD.} \and R. Exel \footnote{ Departamento de
Matem\'atica, UFSC, Florianopolis 88040-900, Brasil. Partially
supported by CNPq, Instituto do Milenio, PRONEX - Sistemas
Din\^amicos.}\and A. O. Lopes\footnote{ Instituto de Matem\'atica,
UFRGS, Porto Alegre 91501-970, Brasil. Partially supported by CNPq,
Instituto do Milenio, PRONEX - Sistemas Din\^amicos, benefici\'ario
de aux\'{\i}lio financeiro CAPES - PROCAD.} } \maketitle

\begin{abstract}

We consider a finite state set $S$ and a continuous time Markov
Chain $X_t$, $t\geq 0$, taking values on $S$. We denote by $\Omega$
the set of paths $w$ taking values on $S$ (the elements $w$ are
locally constant with left and right limits and are also right
continuous on $t$). $P$ will denote the associated probability on
($\Omega, {\cal B}$) which we assume that is stationary. All functions
$f$ we consider bellow are in the set ${\cal L}^\infty (P)$.

From $P$ we are able to define a Ruelle operator ${\cal L}^t,
t\geq 0$, acting on functions $f:\Omega \to \mathbb{R}$ of ${\cal
L}^\infty (P)$. Given $V:\Omega \to \mathbb{R}$, such that is
constant in sets of the form $\{X_0=c\}$, we  define a modified
Ruelle operator ${\cal L}_V^t, t\geq 0$, and we  are able to show
the existence of  an eigenfunction and  an eigen-probability
$\rho_V$ on $\Omega$ associated to ${\cal L}^t_V, t\geq 0$.

We also show the follow property for the  probability $\rho_V$: for
any integrable $f\in {\cal L}^\infty (P)$ and any real and positive $t$

$$
\int e^{ -\int_0^t (V \circ \Theta_s)(.) ds } \, \,[\,( {\cal L}^t
\, (e^{ \int_0^t (V \circ \Theta_s)(.) ds } f )\,)\,   \circ
\theta_t \,] d \rho_V = \int f d \rho_V$$

This equation generalize for continuous time a similar one for
discrete time systems (and which is quite important for
 understanding the KMS states of  certain
$C^*$-algebras).

\end{abstract}

\newpage

\section{Introduction}

We would like to consider a continuous time stochastic process
that maps the positive real line $\RR_+=\{ t \in \RR \colon t \geq 0   \}$
on a finite set $S$ with $n$ elements, that
we can simply write as $S= \{ 1, 2, \ldots , n\}$.
Now take  a $n$ by $n$ real matrix $L$ such that:
\begin{itemize}
\item[1)] $0 < - L_{ii}$, for all $i \in S$,

\item[2)] $ L_{ij} \geq 0,$ for all $i\neq j$, $i \in S$,

\item[3)] $\sum_{i=1}^n L_{ij} =0$ for all fixed
$j \in S.$
\end{itemize}

We point out that, by convention,  we are considering column
stochastic matrices and not line stochastic matrices (see [N]
section 2 and 3 for general references).

We denote by $P^t= e^{t L}$ the semigroup generated by $L$. The
left action of the semigroup can be identified with an action
over functions from $S$ to $\RR$
(vectors in $\mathbb{R}^n$) and the right action
can be identified with action on
measures on $S$ (also vectors in $\mathbb{R}^n$).

The matrix $e^{t L}$ is column stochastic, since from the assumptions on $L$ follows
that
$$
 (1, \ldots, 1) e^{tL} = (1, \ldots, 1)(I+tL+\frac{1}{2}  t^2 L^2  + \cdots) =
  (1, \ldots, 1)
$$

It is well known that there exist a vector of probability $p_0 = (
p_0^1, p_0^2,...,p_0^n ) \in \mathbb{R}^n$ such that $e^{t L}
(p_0)= P^t p_0=p^0$ for all $t>0$. The vector $p_0$ is a right
eigenvector of $e^{t L}$. All entries $p_0^i$ are
{\em strictly positive},  as a consequence of hypothesis 1.

We point out that, if 
$L(p_0)=0$,
then $e^{t L} (p_0)=p_0$, for all $ t\geq 0.$


Now let us consider the space $\tilde{\Om} = \{1,2,...,n\}^{ \RR_+ }$
of all functions from  $\RR_+$ to $S$. In principle
it could be enough for our purposes, but technical details
in the construction of probability measures on such a space
force us to use a restriction of it:
We consider the space $\Omega \subset \tilde{\Om}$ as the
set of right-continuous functions from $\RR_+$ to $S$.
In this set we take the sigma algebra ${\cal B}$ generated by the
cylinders of the form
$$\{ w_{0} = a_0, w_{t_1} = a_1, w_{t_2} = a_2,..., w_{t_r} = a_r\},$$
where $t_i\in \RR_+ , \,r \in \mathbb{Z}^{+},
a_i \in S$ and  $0< t_1 < t_2 < \ldots < t_r$.
It is possible to endow $\Omega$ with a metric,
the Skorohod-Stone metric $d$, which makes
$\Omega$ complete and separable ([EK] section 3.5) but
the space  is not compact.

Now we can introduce a continuous time version of the shift
map as follows:
we define for each fixed $s\in \RR_+$ the
${\cal B}$-measurable transformation $\Theta_s: \Omega \to \Omega$
given by $\Theta_s (w_t) = w_{t+s}$ (we remark  that
$\Te_s$ is also a continuous
transformation with respect to the Skorohod-Stone metric $d$).

For $L$ and $p_0$ fixed as above we denote by $P$ the probability
on the sigma-algebra ${\cal B}$ defined for cylinders by
$$P(\{
w_{0} = a_0, w_{t_1} = a_1,..., w_{t_r} = a_r\})=
   P_{ a_{r} a_{r-1} }^{t_{r}-t_{r-1}}...P_{ a_{2} a_{1} }^{t_2-t_1}
P_{a_{1} a_0 }^{t_1} p_{0}^{a_0}.$$
For details of the construction of this measure the reader is
refered to [B].

The probability $P$ on ($\Omega, {\cal B}$) is stationary in the
sense that for any integrable function  $f$ and any $s \geq 0$
$$ \int f (w)  d P(w)= \int (f \circ \Theta_s) d P (w).$$

From now on the Stationary Process defined by $P$ is
denoted by $X_t$  and
all functions $f$ we consider are in the set ${\cal
L}^\infty (P)$.

 There exist a version of $P$ such that for a set of full
measure all elements $w$ are locally constant on $t$ on the right
side with left and right limits and $w$ is right continuous on
$t$. We consider from now on such $P$.

From $P$ we are able to define a continuous time Ruelle operator
${\cal L}^t$, $t>0$, acting on functions $f:\Omega \to \mathbb{R}$
of ${\cal L}^\infty (P)$. It is also possible to
introduce the endomorphism $\al_t \colon \cL^{\infty}(P)  \to \cL^{\infty}(P)$
defined as
$$
   \al_t(\vfi)=\vfi \circ \Te_t \; , \qquad \forall \vfi \in \cL^{\infty}(P)
$$
Given $V:\Omega \to \mathbb{R}$, such
that it is constant in sets of the form $\{X_0=c\}$
(i.e., $V$ depends only on the value of $x(0)$), we are able to
show the existence of a probability $\rho_V$ on $\Omega$ which
is absolutely continuous with respect to $P$ and satisfies:

\begin{theoA}
 For any integrable $f\in {\cal L}^\infty (P)$ and any positive
$t$
$$
\int e^{ -\int_0^t (V \circ \Theta_s)(.) ds } \, \,[\,( {\cal L}^t
\, (e^{ \int_0^t (V \circ \Theta_s)(.) ds } f )\,)\,   \circ
\theta_t \,] d \rho_V = \int f d \rho_V
$$
\end{theoA}

The above functional equation is a natural generalization (for
continuous time) of the similar one presented in [EL1] and [EL2]. We
believe it will be important in the analysis of certain $C^*$
algebras, generated by the operators $\al$ and $\cL$, specially concerning
the characterization of KMS states. We point out however that we are able to show this
property of $\rho_V$
just for a quite simple function $V$ as above.

With the operators $\al$ and $\cL$  we can rewrite the theorem above as
$$
  \rho_V(G_T^{-1} E_T (G_T \vfi) )= \rho_V(\vfi)
$$
for all $\vfi \in \cL^{\infty}$ and all $T > 0$,
where, as usual,  $\rho_V(\vfi)=\int \vfi d\rho_V$,
$E_T = \al_T \cL^T$ is in fact a projection on a subalgebra of $\cB$
and
$G_T \colon \Om \to \RR$ is given by
$$
 G_T(x)= \exp{(\int_0^T V(x(s))ds )}
$$

For the map $V : \Omega \to \mathbb{R}$, which is constant in
cylinders of the form $\{ w_{0} = i\}$, $i\in \{1,2,...,n\}$, we
denote by $V_i$ the corresponding value. We denote also by $V$ the
diagonal matrix with the $i$-diagonal element equal to $V_i$.

We denote by $P_V^t= e^{t (L +V)}.$ The Perron-Frobenius Theorem
for such semigroup will be one of the main ingredients of the
proof.

A related and more general result will appear in [LNT].

As usual we denote by ${\cal F}_s$ the sigma-algebra generated by
$X_s$. We also denote by ${\cal F}_s^{+}$ the sigma-algebra
generated $\sigma( \{ X_u, s\leq u   \})$. Note that a ${\cal
F}_s^{+}$-measurable function $f(w)$ on $\Omega$ does  depend of
the value $w_s$.

We also denote by $I_A$ the indicator function of a measurable set
$A$ in $\Omega$.

\vspace{0.3cm}

\section{A continuous time  Ruelle Operator}


We consider the disintegration of $P$ given by the
 family of measures, indexed by the elements of $\Om$ and $t > 0$
defined as follows: first,
consider a sequence $0 = t_0 < t_1 < \cdots < t_{j-1} < t \leq t_j
< \cdots < t_r$. Then for $w \in \Om$ and $t > 0$
 we have   on cylinders:
$$
  \mu^w_t( [ X_0=a_0, \ldots, X_{t_r}=a_r]) =
$$
$$
 \left\{
   \begin{array}{cc}
    \frac{1}{p_0^{w(t)}} P_{w(t) a_{j-1}}^{t-t_{j-1}} \cdots
         P_{a_2 a_1}^{t_2-t_1} P_{a_1 a_0}^t p_0^{a_0}  &
    \text{if $a_j=w(t_j), \ldots, a_r=w(t_r)$}  \\
     0 & \text{otherwise}
   \end{array}
 \right.
$$

\begin{prop}
    $\mu^w_t$ is the disintegration of $P$.
\end{prop}
\begin{demo}
  It is enough to show that for any integrable $f$
$$
   \int_{\Om} f dP = \int_{\Om}
   \int_{\Te_t^{-1}(w)} f(x) d\mu^w_t(x) dP(w)
$$
For doing that we can assume that $f$ is in fact the indicator
of the cylinder $[X_0=a_0, \ldots, X_{t_r}=a_r]$; then the
right hand side becomes
$$
  \int \int f d\mu^w_t(x) dP(w)= \int I_{[w(t_j)=a_j, \ldots, w(t_r)=a_r ]}
   \frac{1}{p_0^{w(t)}} P^{t-t_{j-1}}_{w(t) a_{j-1}} \cdots
   P_{a_1 a_0}^{t_1} p_0^{a_0}  dP(w) =
$$
$$
  \sum_{a=1}^n \int I_{[w(t)=a, w(t_j)=a_j, \ldots, w(t_r)=a_r ]}
 \frac{1}{p_0^{w(t)}} P^{t-t_{j-1}}_{w(t) a_{j-1}} \cdots
   P_{a_1 a_0}^{t_1} p_0^{a_0}  dP(w) =
$$
$$
\sum_{a=1}^n P^{t_r-t_{r-1}}_{a_r a_{r-1}} \ldots P^{t_j-t}_{a_{t_j} a} p_0^a
  \frac{1}{p_0^a} P^{t-t_{j-1}}_{a a_{j-1}} \cdots
   P_{a_1 a_0}^{t_1} p_0^{a_0} =
P^{t_r-t_{r-1}}_{a_r a_{r-1}}  \cdots
   P_{a_1 a_0}^{t_1} p_0^{a_0} =
$$
$$
  P([X_0=a_0, \ldots, X_{t_r} =a_r]) = \int f dP
$$
The proof for a general $f$ now follows
from standard arguments.
\end{demo}

\begin{defi}
For $t$ fixed we define the operator $\cL^t: \cL^\infty
(\Omega,P) \to \cL^\infty (\Omega,P) $ as follows:
$$
   \cL^t(\vfi)(x)= \int_{\bar{y} \in \Te_t^{-1}(x)}
            \vfi(\bar{y}) d\mu^x_t(\bar{y})
$$
\end{defi}

\begin{rema}
The definition above can be rewritten as
$$
\cL^t(\vfi)(x) =  \int_{y \in D[0, t) } \vfi(yx) d\mu^x_t(yx)
$$
where the symbol $yx$ means the concatenation of the
path $y$ with the translation of $x$:
$$
     xy(s)=
  \left\{
    \begin{array}{cc}
       y(s)  &  \text{if $s \in [0, t)$}    \\
        x(s-t) &  \text{if $s \geq  t$}
    \end{array}
  \right.
$$
and $D[0, t)$ is the set of right-continuous functions
from $[0, t)$ to $S$. This follows simply from the fact
that, in this notation, $\Te_t^{-1}(x)= \{  yx \colon y \in D[0, t) \}$.
\end{rema}


It is possible to shed some light on the meaning of this
operator applying it to some simple functions.
For example,  we can  see the effect
of $\cL^t$  on some indicator of a given cylinder:
Consider the sequence  $0=t_0< t_1<..< t_{j-1}< t \leq t_j< ...
<t_r $ and then
 take
 $f= I_{ \{ X_{0} = a_0, X_{t_1} = a_1,..., X_{t_r}
= a_r\}}$. Then, for a path $z \in \Om$  such that $z_{t_j-t}=a_j,
..., z_{t_r-t} = a_r$ (the future condition)
we have
$$  {\cal L}^t (f) (z) =\frac{1}{p^{z_0}_0} P_{ z_{0} a_{j-1}
}^{t -t_{j-1}}...P_{ a_{2} a_{1} }^{t_2-t_1} P_{a_{1} a_0 }^{t_1}
p_{0}^{a_0},$$ otherwise (i.e., if the path $z$ does not satisfy the
condition above) we get ${\cal L}^t (f) (z)=0.$


Note that if $t_r < t$, then ${\cal L}^t (f) (z)$ depends only on
$z_0$.   For example, if $ f=I_{\{X_0=i_0\}}$ then
$$
 {\cal L}^t (f) (z)= \int_{y \in D[0, t)} I_{\{X_0=i_0\}}(yx) d\mu^{z}_t(yx)=
 \mu^z_t([X_0=i_0])=
 \frac{1}{p_0^{z_0} }   P^t_{z_0, i_0} p_0^{i_0}
$$
In the case  $
f=I_{\{X_0=i_0, X_t=j_0\}}$, then ${\cal L}^t (f) (z)= P^t_{z_0,
i_0} \frac{p_0^{i_0}}{p_0^{z_0} }$, if $z_0=j_0$, and ${\cal L}^t
(f) (z)=0$ otherwise.

Now we can show some properties of $\cL^t$.

\begin{prop}
 ${\cal L}^t (1)=1$, where
$1$ is the function that maps every point in $\Om$ to $1$.
\end{prop}
\begin{demo}
  Indeed
$$
  \cL^t(1)(x)= \int_{y \in D[0, t)} 1(yx) d\mu^x_t(yx) =
    \int d\mu^x_t(yx) = \mu^x_t([X_t = x(0)])=
$$
$$
   \sum_{a=1}^n \mu^x_t([X_0 = a, X_t = x(0)]) =
 \frac{1}{p_0^{x(0)}} \sum_{a = 1}^n P^t_{x(0)  a} p_0^a = 1
$$
\end{demo}

We can  also define the dual of $\cL^t$, denoted by $(\cL^t)^*$,
acting on the measures. Then we get:



\begin{prop}
For any positive $t$ we have that $({\cal L}^t)^* (P) = (P)$
\end{prop}

\begin{demo} For a fixed $t$ we have that $({\cal L}^t)^* (P) = (P)$
because for any  $f$ of the form
$f= I_{ \{ X_{0} = a_0, X_{t_1} =
a_1,..., X_{t_r} = a_r\}}$, $0=t_0< t_1<..< t_{j-1}< t \leq t_j<
... <t_r.  $
 we have
$$
\int {\cal L}^t (f)(z)  d P(z) = \sum_{b=1}^n \int_{\{X_0=b\} }
{\cal L}^t (f)(z)  d P(z)=
$$
$$
 \sum_{b=1}^n \int I_{\{X_0= b,
X_{t_j-t}=a_j, ..., X_{t_r-t} = a_r \} } (z) d P(z)
\frac{1}{p^{b}_0} P_{ b a_{j-1} }^{t -t_{j-1}}...P_{ a_{2} a_{1}
}^{t_2-t_1} P_{a_{1} a_0 }^{t_1} p_{0}^{a_0}=
$$
$$
\sum_{b=1}^n P(\{X_0= b,X_{t_j-t}=a_j, ..., X_{t_r-t} =
a_r \})\frac{1}{p^{b}_0} P_{ b a_{j-1} }^{t
-t_{j-1}}...P_{ a_{2} a_{1} }^{t_2-t_1} P_{a_{1} a_0 }^{t_1}
p_{0}^{a_0}=
$$
$$
 \int f (w) d P(w).
$$
\end{demo}

\begin{prop}
 Given  $t \in \RR_+$ and the functions  $\vfi, \psi \in \cL^{\infty}(P)$
then  we have
$$
{\cal L}^t (\vfi \times (\psi \circ
\Theta_t) ) (z)= \psi(z) \times {\cal L}^t (\vfi)(z).
$$
\end{prop}

\begin{demo}
$$
   \cL^t(\vfi (\psi \circ \Te_t))(x)=
    \int_{i \in D[0, t)} \vfi(ix) (\psi \circ \Te_t)(ix) d\mu^x_t(i) =
$$
$$
 \psi(x)  \int \vfi(ix) d\mu^x_t(i) = (\psi \cL^t(\vfi))(x) =
\psi(x) \cL^t(\vfi)(x)
$$
since $\psi \circ \Te_t(ix) = \psi(x)$, independently of $i$.
\end{demo}

We just recall that the last proposition can be restated as
$$
     \cL^t(\vfi \al_t(\psi) )= \psi \cL^t(\vfi)
$$

Then we get:
\begin{prop}
 $\al_t $ is the dual of $ {\cal L}^t$ on ${\cal L}^2 (P)$.
\end{prop}

\begin{demo}
 From last two propositions
$$ \int  {\cal L}^t (f) g   \,d P= \int {\cal L}^t
( f \times ( g \circ \Theta_t))\, d P =
\int f \times ( g \circ \Theta_t)\, dP = \int f \al_t(g) dP  $$
as claimed.
\end{demo}

Now we would like to obtain conditional expectations.
 For a given
$f$ recall that the function $Z (w)= E ( f  |{\cal F}_t^{+})$ is
the $Z$   (almost everywhere defined) ${\cal F}_t^{+}$-measurable
function  such that for any ${\cal F}_t^{+}$-measurable set $B$ we
have $\int_B Z(w) d P(w) = \int_B f(w) d P(w).$

\begin{prop}
  The conditional expectation is given by
$$
   E(f | \cF_t^+)(x) = \int f d\mu_t^x
$$
\end{prop}
\begin{demo}
 For $t$ fixed, consider a  ${\cal F}_t^{+}$-measurable set $B$. Then  we
have

$$ \int_B E(f| \cF_t^+) dP =
 \int_B \int f d \mu_t^w d P(w)= \int (I_B (w) \int f d \mu_t^w) d
 P(w)=
$$
$$ \int  \int (f I_B) d \mu_t^w d P(w)=  \int f (w) I_B (w)
d P(w)= \int f dP,
$$
and the proposition is concluded.
\end{demo}

Now we can relate the conditional expectation with respect to the
$\si$-algebras $\cF_t^+$ with the operators $cL^t$ and $\al_t$
as follows:
\begin{prop}
 $[{\cal L}^t (f) ]( \Theta_t ) = E ( f |{\cal F}_t^{+})$ (i.e. $E=\al_t  \cL^t$).
\end{prop}
\begin{demo}
This follows from the fact that for any $B=\{ X_{s_1} = b_1,
X_{s_2} = b_2,..., X_{s_u} = b_u\},$ with $t< s_1<... <s_u,$ we
have $I_B = I_A \circ \Theta_t$ for some measurable $A$ and
$$ \int_B {\cal L}^t (f ) (\Theta_t  (w)) d P(w)=
\int I_B (w) {\cal L}^t (f) ( \Theta_t  (w)) d P(w)= $$
$$\int (I_A
\circ \Theta_t)(w) {\cal L}^t (f ) ( \Theta_t  (w)) d P(w)=
 \int I_A (w) {\cal L}^t (f )(w)  d P(w)$$
$$\int {\cal L}^t (f (I_A \circ \Theta_t) ) (w) d P(w)= \int f (w)I_A  (\Theta_t (w))
d P(w)= \int_B f(w) d P(w)$$
\end{demo}

\section{The modified operator}

We are interested in the perturbation by $V$ (defined above) of
the ${\cal L}^t$ operator.

\begin{defi}
  We define $G_t \colon \Om \to \RR$ as
$$
    G_t(x)= \exp{(\int_0^t V(x(s))ds  )}
$$
\end{defi}

\begin{defi}
We define the $G$-weigthed transfer operator
${\cal L}^t_V:  {\cal L}^\infty (\Omega,P) \to {\cal
L}^\infty (\Omega,P)  $ acting on measurable functions $f$ (of the
above form) by
$$
  \cL_V^t(f)(w):= \cL^t(G_t f) =
$$
$$={\cal L}^t (  e^{ \int_0^t (V \circ \Theta_s)(.) ds } \,  f \, ) =
\sum_{b=1}^n {\cal L}^t (  e^{ \int_0^t (V \circ
\Theta_s)(.) ds } \, I_{\{X_t =b\}} \, f \, ) (w)
$$
\end{defi}

Note that $e^{ \int_0^t (V \circ \Theta_s)(.) ds } \, I_{\{X_t
=b\}} \, $ does not depend on information larger then $t$. In the
case $f$ is such that $t_r\leq t$ (in the above notation), then
${\cal L}^t_V (f) (w)$ depends only on $w(0)$.

The integration on $s$ above is over the open interval $(0,t)$.

Remember that if
$L(p_0)=0$,
then $e^{t L} (p_0)=p_0$ for all $ t\geq 0.$

We will consider soon  an eigenfunction and an eigen-measure for
the operator ${\cal L}^t_V$. But, first we need the following:
\vspace{0.2cm}

\begin{theo}
   ([S] page 111) We assume $S$ is finite. One can prove that for $L, p_0$ (such that $L(p_0)=0$), and $V$ fixed as above there exists
\begin{itemize}
\item[a)] a unique strictly positive function $u_V : \Omega\to \mathbb{R}$,
constant equal  to the value $u_V^i$ in each cylinder $X_0= i$, $i
\in \{1,2,..,n\},$ (we can see $u_V$ as  $u_V:S \to \mathbb{R}$, or,  as a vector
in $\mathbb{R}^n$),

\item[b)]  a unique probability vector  $\mu_V$ in $\mathbb{R}^n$(a
probability over over the set $ \{1,2,..,n\}$ such that
$\mu_V(\{i\})>0,\, \forall i$), such that
$$ \sum_{i=1}^n
u_V^i (\mu_V)_i = 1 ,$$

\item[c)] a real positive value $\lambda (V)$,

such that

\item[d)] for any positive $s$
$$ e^{-s \lambda (V)} u_V \, e^{s (L+ V)} \, = u_V. $$

Moreover, for any  $v =(v_1,.,v_n) \in \mathbb{R}^n$
$$\lim_{t\to \infty}  e^{-t \lambda (V)} \, v \,  e^{t (L + V)}  =
(\sum_{i=1}^n v_i (\mu_V)_i) \, u_V, $$

\item[e)] for any positive $t$
$$ (P^t_{V}) \, \mu_V = e^{\lambda(V)\,t} \mu_V. $$
\end{itemize}
\end{theo}

From property e) it follows (a right eigenvalue) that
$$(L + V) \,\mu_V  = \lambda (V) \, \mu_V.$$
From d) it follows (a left eigenvalue) that
$$u_V\, (L + V)   = \lambda (V) \, u_V.$$

 Note that when $V=0$, then $\lambda(V)=0$, $\mu_V=p_0$ and $u_V$
 is constant equal to $1$.

\vspace{0.2cm}

{\bf In order to show the existence of $u_V$, such that, $u_V (L + V) \, = \lambda (V) \, u_V$ one add a constant to $V$ in such way that all the entries of $(L + V)$ are positive. By the Perron Theorem, this will imply d), that is, the existence of $\lambda (V) \in \mathbb{R}$ and a vector $u_V$, such that,  $u_V\, (L + V)   = \lambda (V) \, u_V.$

\bigskip

For the case the space $S$ is not finite see [LNT].}
\bigskip

Now we return to our setting: for each $i_0$ and $t$ fixed one can
consider the probability $\mu_{i_0}^t$ defined over the
sigma-algebra ${\cal F}_t^{-}= \sigma ( \{ X_s| s \leq t \} )$
with support on $\{X_0=i_0\}$ such that for cylinder sets with $0<
t_1<...<t_r\leq t$
$$\mu_{i_0}^t(\{
X_{0} = i_0, X_{t_1} = a_1,..., X_{t_{r-1}} = a_{r-1}, X_t =
j_0\})= P_{ j_0 a_{r} }^{t-t_{r}}...P_{ a_{2} a_{1} }^{t_2-t_1}
P_{a_{1} i_0 }^{t_1} .$$

The probability $\mu_{i_0}^t$ is not stationary.

We denote by $Q(j, i)_t$ the $i,j$ entry of the matrix $e^{t (L +
V)}$, that is $(e^{t (L + V)})_{i,j}$.

It is known ([K] page 52 or  [S] Lemma 5.15) that
$$ Q(j_0, i_0)_t= E_{\{X_0= i_0\}}
\{ e^{ \int_0^t (V \circ \Theta_s)(w) ds } \,;\, X(t) = j_0\} =$$
$$
\int \, I_{\{X_t = j_0\}}\,\, e^{ \int_0^t\,(V \circ \Theta_s)(w)
ds } d \mu_{i_0}^t (w).$$

For example,
$$ \int I_{\{X_t = j_0\}} e^{ \int_0^t (V \circ \Theta_s)(w) ds } d P=
\sum_{i=1,2,..,n } Q(j_0, i)_t\,  p_i^0$$

In the particular case where  $V$ is constant equal $0$,
then $p^0 = \mu_V$ and
$\lambda(V)=0$.

\begin{prop}
$f(w)=\frac{\mu_V(w)}{ p^0 (w)}= \frac{(\mu_{V})_{ w(0)}}{
(p^0)_{w(0)}  } $ is an eigenfunction for ${\cal L}^t_V$
with eigenvalue $e^{t\la(V)}$.
\end{prop}

\begin{demo} Note that $\frac{\mu_V}{ p^0 }=\sum_{c=1}^n \frac{\mu_V(c)}{
p^0 (c) }\, I_{\{ X_0 = c\} }.$

For a given $w$, denote $w(0)$ by $j_0$, then conditioning
$${\cal L}^t_V (\frac{\mu_V}{ p^0 }) (w)= \sum_{ c=1}^n \,\sum_{b=1}^n
{\cal L}^t_V\, ( \, \frac{\mu_V(c)}{ p^0 (c) }\, I_{\{ X_0 = c\}
}\,I_{\{ X_t = b\} }\,)\, (w).
$$

Consider $c$ fixed, then for $b=j_0$ we have $${\cal L}^t_V\, ( \,
I_{\{ X_0 = c\} }\,I_{\{ X_t = b\} }\,)\, (w)= \frac{Q(j_0,c)_t \,
\, p^0_c }{p^0_{j_0}},
$$
and for $b\neq j_0$, we have ${\cal L}^t_V\, ( \, I_{\{ X_0 = c\}
}\,I_{\{ X_t = b\} }\,)\, (w)=0$.

Finally,
$$ {\cal L}^t_V (\frac{\mu_V}{ p^0 }) (w)= \sum_{c=1}^{n} \frac{\mu_V(c)}{ p^0 (c) }
\,  Q(j_0, c)_t \, \, \frac{ p^0 _{c}}{ p^0 _{j_0}}=
 e^{t \lambda(V)} \frac{(\mu_V)_{j_0}}{ p^0_{j_0}}= e^{t \lambda(V)} (\frac{\mu_ V}{
p^0})(w),
$$
because $ \, e^{t (L + V)}\, (\mu_V) = e^{t \lambda(V)} (\mu_V)$.

Therefore for any $t>0$ the function $\frac{\mu_V}{ p^0 }$ (that
depends only on $w(0)$) is an eigenfunction for the operator
${\cal L}^t_V $ associated to the eigenvector $ e^{t
\lambda(V)}$.
\end{demo}

\vspace{0.2cm}

\begin{defi}
Consider now for each $t$ the operator acting on $g$ by
$$\hat{{\cal L}}^t_V  (g) (w)= [\, \frac{p^0 }{\mu_V  }
{\cal L}^t (e^{ \int_0^t (V -  \lambda (V)) \circ \Theta_s)(.) ds
} \, g \, \frac{\mu_V}{ p^0 } ) \,]\, (w)
$$
\end{defi}

From the above $\hat{{\cal L}}^t_V  (1)=1$ for all positive $t$.

Note that by conditioning, if
$g=I_{\{X_0=a_0,X_{t_1}=a_1,X_{t_2}=a_2,X_{t}=a_3\}}$, with
$0<t_1< t_2< t$, then
$$ \hat{{\cal L}}^t_V (g) (w)= \frac{\mu_V (a_0) }{\mu_V ( a_3)}
e^{(t-t_2)\, (L+V-\lambda\, I)}_{a_3 a_2 }\, e^{(t_2-t_1)\,
(L+V-\lambda\, I)}_{a_2 a_1 }\, e^{t_1\,(L+V-\lambda\, I)}_{a_1
a_0 },$$ for $w$  such that $w_0=a_3$, and $\hat{{\cal L}}^t_V (g)
(w)=0$ otherwise.

Moreover, for
$g=I_{\{X_0=a_0,X_{t_1}=a_1,X_{t}=a_2,X_{t_3}=a_3\}}$, with
$0<t_1< t< t_3$, then
$$ \hat{{\cal L}}^t_V (g) (w)= \frac{\mu_V (a_0) }{\mu_V ( a_2)}
e^{(t-t_1)\, (L+V-\lambda\, I)}_{a_2 a_1 }\,
e^{t_1\,(L+V-\lambda\, I)}_{a_1 a_0 },$$ for $w$  such that
$w_0=a_2, w_{t_3-t}= a_3$, and $\hat{{\cal L}}^t_V (g) (w)=0$
otherwise.

Consider now the dual operator $(\hat{{\cal L}}^t_V)^*$.

For $t$ fixed consider the transformation in the set of
probabilities $\mu$ on $\Omega$ given by $(\hat{{\cal L}}^t_V)^*
(\mu) = \nu$.

\begin{theo}
  There exists a fixed probability measure $\nu_V$ on
$(\Omega,{\cal B})$ for such transformation $(\hat{{\cal
L}}^t_V)^*$. The stationary probability $\nu_V$ does not depend on
$t$.
\end{theo}
\begin{demo}
Denote by $\nu=\nu_V$ the probability obtained in the following
way, for
$$g=I_{\{X_0=a_0,X_{t_1}=a_1,X_{t_2}=a_2,\,...,X_{t_{r-1}}=a_{r-1},X_{r}=a_r\}},$$
with $0<t_1<t_2<...<t_{s-1}< t\leq t_{s}<..<t_r$, we define
$$ \int g(w) \, d \nu(w)=
e^{(t_r-t_{r-1})\, (L+V-\lambda\, I)}_{a_r a_{r-1}}\,...
\,e^{(t_2-t_1)\, (L+V-\lambda\, I)}_{a_2 a_1 }\,
e^{t_1\,(L+V-\lambda\,I)}_{a_1 a_0 }\, \mu_V (a_0).$$

This probability satisfies the Kolmogorov compatibility conditions
because is defined via a semigroup (see chapter IV. 2 [BW])

In order to show that $\nu$ is a probability  we have to  use the
fact that $\sum_{c\in \tilde{S}}\, \mu_V (c)=1$

On the other hand,
$$z(w)= \hat{{\cal L}}^t_V (g) (w)= \frac{\mu_V (a_0) }{\mu_V ( w_0)}
e^{(t-t_{s-1})\, (L+V-\lambda\, I)}_{w_0 a_{s-1}}\,...
\,e^{(t_2-t_1)\, (L+V-\lambda\, I)}_{a_2 a_1 }\,
e^{t_1\,(L+V-\lambda\,I)}_{a_1 a_0 },$$ for $w$  such that $w_{t_s
- t}=a_s,w_{t_{s+1} - t}=a_{s+1},...,w_{t_r - t}=a_r$, and
$\hat{{\cal L}}^t_V (g) (w)=0$ otherwise. Note that $
z(w)=\hat{{\cal L}}^t_V (g) (w)$ depends only on $w_0,w_{t_{s+1} -
t},...,w_{t_r - t}$.

We have to show that for any $g$ we have $\int g\, d \nu= \int
\hat{{\cal L}}^t_V (g) d \nu$.

 Now,
$$ \int z(w) \, d \nu (w)=\int \sum_{c\in S} I_{\{X_0=c,X_{t_s
- t}=a_s,X_{t_{s+1} - t}=a_{s+1},...,X_{t_r - t}=a_r   \}}\,z(w)
\,\, d \nu (w)= $$
$$\sum_{c\in S}\nu (\{X_0=c,X_{t_s
- t}=a_s,X_{t_{s+1} - t}=a_{s+1},...,X_{t_r - t}=a_r \})$$
$$\frac{\mu_V (a_0) }{\mu_V ( c)} e^{(t-t_{s-1})\, (L+V-\lambda\,
I)}_{c a_{s-1}}\,... \,e^{(t_2-t_1)\, (L+V-\lambda\, I)}_{a_2 a_1
}\, e^{t_1\,(L+V-\lambda\,I)}_{a_1 a_0 }= $$
$$\sum_{c\in S} \mu_V (c)  \,
e^{(t_r-t_{r-1})\, (L+V-\lambda\, I)}_{a_r a_{r-1}}\,...
\,e^{(t_{s+1} -t_s)\, (L+V-\lambda\, I)}_{a_{s+1} a_s }\,
e^{t_s-t\,(L+V-\lambda\,I)}_{a_{s} c }$$
$$\frac{\mu_V (a_0)}{\mu_V (c)}  \,
e^{(t-t_{s-1})\, (L+V-\lambda\, I)}_{c a_{s-1}}\,...
\,e^{(t_2-t_1)\, (L+V-\lambda\, I)}_{a_2 a_1 }\,
e^{t_1\,(L+V-\lambda\,I)}_{a_1 a_0 }=$$
$$\sum_{c\in S}   \,
e^{(t_r-t_{r-1})\, (L+V-\lambda\, I)}_{a_r a_{r-1}}\,...
\,e^{(t_{s+1} -t_s)\, (L+V-\lambda\, I)}_{a_{s+1} a_s }\,
e^{t_s-t\,(L+V-\lambda\,I)}_{a_{s} c }$$
$$\mu_V (a_0)  \,
e^{(t-t_{s-1})\, (L+V-\lambda\, I)}_{c a_{s-1}}\,...
\,e^{(t_2-t_1)\, (L+V-\lambda\, I)}_{a_2 a_1 }\,
e^{t_1\,(L+V-\lambda\,I)}_{a_1 a_0 }=$$
$$\mu_V (a_0)\, e^{(t_r-t_{r-1})\, (L+V-\lambda\, I)}_{a_r a_{r-1}}\,...
\,e^{(t_{s+1} -t_s)\, (L+V-\lambda\, I)}_{a_{s+1} a_s }\,$$
$$( \sum_{c\in S}\,e^{(t_s-t)\,(L+V-\lambda\,I)}_{a_{s} c }
e^{(t-t_{s-1})\, (L+V-\lambda\, I)}_{c a_{s-1}})\, \,...
\,e^{(t_2-t_1)\, (L+V-\lambda\, I)}_{a_2 a_1 }\,
e^{t_1\,(L+V-\lambda\,I)}_{a_1 a_0 }=
$$
$$\mu_V (a_0)\, e^{(t_r-t_{r-1})\, (L+V-\lambda\, I)}_{a_r a_{r-1}}\,...
\,e^{(t_{s+1} -t_s)\, (L+V-\lambda\, I)}_{a_{s+1} a_s }\,$$
$$e^{(t_s-t_{s-1})\,(L+V-\lambda\,I)}_{a_{s} a_{s-1} }\, \,... \,e^{(t_2-t_1)\, (L+V-\lambda\, I)}_{a_2 a_1 }\,
e^{t_1\,(L+V-\lambda\,I)}_{a_1 a_0 }=
$$

$$\int g\, d \nu .$$

The claim for the general $g$ follows from the above result.
\end{demo}

 \vspace{0.3cm}

\begin{defi}
Consider  the probability $\rho_V= (g_V)^{-1} \nu_V$, where $g_V$
is chosen colinear to $\frac{\mu_V}{ p^0 }$, in such way $\rho_V$
is a probability (not necessarily invariant) on $\Omega$.
\end{defi}

It easily follows that $({\cal L}^t_V)^* (\rho_V) = e^{t
\lambda_V} \rho_V.$ The probability $\nu_V$ is invariant for
$\theta_s$ with $s\geq 0$.

From last theorem  follows easily: \vspace{0.3cm}

\begin{prop}
 For any integrable $f, g \in {\cal L}^\infty (P)$ and any positive $t$
$$
  \int g {\cal L}^t_V (f) d \rho_V = \int {\cal L}^t_V (f (g \circ \theta_t) ) d \rho_V
 =  e^{t \lambda_V}  \int f (g \circ \theta_t)  d \rho_V.
$$
\end{prop}

Now we are in position to prove our main result.
From $({\cal L}^t_V)^* (\rho_V) = e^{t \lambda_V} \rho_V$ it
follows that the measure $\rho_V$ satisfies the important
equation: \vspace{0.3cm}

\begin{theoA}
For any integrable $f\in {\cal L}^\infty (P)$ and any positive $t$
$$
\int e^{ -\int_0^t (V \circ \Theta_s)(.) ds } \, \,[\,( {\cal L}^t
\, (e^{ \int_0^t (V \circ \Theta_s)(.) ds } f )\,)\,   \circ
\theta_t \,] d \rho_V = \int f d \rho_V$$
\end{theoA}

\begin{demo}
$$
\int e^{ -\int_0^t (V \circ \Theta_s)(.) ds } \, \,[\,( {\cal L}^t
\, (e^{ \int_0^t (V \circ \Theta_s)(.) ds } f )\,)\,   \circ
\theta_t \,] d \rho_V = $$
$$e^{-t \lambda_V}\, \int [{\cal L}^t
(e^{ -\int_0^t (V \circ \Theta_s)(.) ds } e^{ \int_0^t (V \circ
\Theta_s)(.) ds })] \, [\, {\cal L}^t \, (e^{ \int_0^t (V \circ
\Theta_s)(.) ds } f )\,]\, d \rho_V=
$$
$$ e^{-t \lambda_V}\, \int  {\cal L}^t
\, (e^{ \int_0^t (V \circ \Theta_s)(.) ds } f ) d \rho_V=
$$
$$ \int f d \rho_V
$$
\end{demo}

\vspace{1.0cm}



\bigskip

\centerline{\bf Bibliography} \vspace{1.0cm}

[B] P. Billingsley, Convergence of probability measures, John Wiley, 1968

[BW] R. Bhattacharya and E. Waymire, Stochastic Processes with
applications, Wiley (1990)

[EK] S. Ethier and T. Kurtz, Markov Processes, John Wiley, 1986

[EL1] R. Exel and A. O. Lopes, "$C^{*}$-Algebras, approximately
proper equivalence relations, and Thermodynamic Formalism", {\it
Erg. Theo. and Dyn. Systems.} {\bf  24} (2004) 1051-1082

[EL2] R. Exel and A. O. Lopes, ''$C^{*}$-Algebras and Thermodynamic
Formalism",  {\it S\~ao Paulo Journal of  Mathematical Sciences}  (USP) 2, 1 (2008), 285–-307

[K] M. Kac, Integration in Function spaces and some of its
applications, Acad Naz dei Lincei Scuola Superiore Normale
Superiore, Piza, Italy (1980).

[LNT] A. Lopes, A. Neumann and P. Thieullen, A Thermodynamic Formalism for continuous time Markov chains with values on the Bernoulli Space: entropy, pressure and large deviations, to appear

[N] J. B. Norris, Markov Chains,  Cambridge Press

[P] K. Parthasarathy, Probability measures on metric spaces,
Academic Press,

[S] D. W. Strook, An introduction to the Theory of Large Deviations, Springer  ,
1984

\end{document}